# BENFORD'S LAW FROM 1881 TO 2006:

# A BIBLIOGRAPHY

## EDITED BY

## WERNER HÜRLIMANN

Feldstrasse 145
CH-8004 Zürich
E-mail : whurlimann@bluewin.ch
Homepage : www.geocities.com/hurlimann53

July 5, 2006